\documentclass[12pt]{article}
\usepackage{amsmath}
\usepackage{amsfonts}
\usepackage{graphicx}

\title{Analytic matrix technique for boundary value problems in applied plasticity  }

\author{L. Novozhilova\\
 {\it Department of Mathematics,}\\
  {\it Western Connecticut State University,
   Danbury, CT, USA\/}\\
S. Urazhdin\\
{\it Department of Physics and Astronomy,}\\
 {\it West Virginia University,
  Morgantown, WV, USA\/}
}
\date{}
\begin{document}
\maketitle
\begin{abstract}
\noindent

An efficient matrix formalism  for finding  power series solutions
to boundary value problems  typical for technological plasticity
is developed. Hyperbolic system of two first order quasilinear
PDEs  that models two-dimensional plastic flow  of von Mises
material  is converted to  the telegraph equation by the hodograph
transformation. Solutions to the  boundary value problems are
found in terms of hypergeometric functions. Convergence issue is
also addressed. The method is illustrated by two test problems of
metal forming.

\end{abstract}

\section{Introduction}

A model of incompressible, ideal, rigid-plastic material was
developed about a hundred years ago. In the two-dimensional case
the model is described by a hyperbolic system of quasilinear
equations. In 1920s Prandl, Hencky, and Mises, among others,
suggested a systematic way, named slip-line method (SLM), for
finding stress fields and associated velocity fields for
two-dimensional plain strain deformation problems for
rigid-plastic body . The method is presented in many classic books
on plasticity theory (cf. \cite{hill}, \cite{kachanov},
\cite{lubliner}) along with its applications to metal forming
processes and problems of plastic failure (limit analysis).
Solution of a 2D rigid-plastic flow problem by this method is
reduced to construction of the field of characteristics
(slip-lines). Assumed domain of the plastic flow is to be
decomposed into a set of subdomains (patches) of a priori
specified types, and appropriate boundary value problems (BVPs),
consistent with physics of the process, are stated for each of the
subdomains. The subproblems were usually solved numerically. In
1980s a matrix implementation of the SLM, called matrix method,
was developed within engineering community \cite{johnson}. The
matrix method is based on the assumption that boundary data are
defined by real analytic functions. The solution to each
subproblem is obtained by applying appropriate matrix operators to
the vector(s) of coefficients of the data. In \cite{novo} a new,
simpler approach to implementation of the matrix method was
developed, mathematically justified, and equipped with recurrent
definitions of five matrix operators needed for solving main
boundary value problems typical in technological applications
(versus twenty operators in the original engineering version of
the method). It was also shown that solutions to the subproblems
can be written in terms of hypergeometric functions.

In this work   explicit expressions for the five matrix operators
of the slip-line analytic technique (SLAT) developed in
\cite{novo} are presented. Convergence issues are addressed and
model examples are provided.

SLAT can be used as a source of test problems for numerical
methods aimed at problems with discontinuities. The fact that
solutions to problems in plasticity may have singularities along
characteristics is well known and experimentally confirmed.
Presence of discontinuities is an important feature  that cannot
be easily detected by numerical methods. In general, capturing
singularities is very challenging mathematical problem that is of
great interest in many applications, including phase transitions
and microstructure formation. SLAT also provides examples of
analytically solved problems with free boundaries since only a
structure of the plastic domain is "guessed" a priory and finding
the boundary of the domain is part of the solution. Furthermore,
although analytic solutions have been superseded with powerful
numerical methods, in particular FEM, the slip-line analysis is a
working tool in metal forming \cite{nepershin}, granular flow
modelling \cite{drescher}, and geomechanics \cite{davis}. The
velocity field in the hodograph plane can also be found by this
method \cite{novo}. Exact solutions for stress and velocity fields
can be used to derive analytic expressions for the plastic power,
which provides a framework for finding optimal in a certain sense
geometric or material parameters of the process under
consideration.

SLM is also a classic topic in standard courses on plasticity
theory, and so is the method of characteristics in PDE theory.
Therefore the authors believe that an elegant matrix formulation
of the method presented in this paper has a good potential in
engineering and mathematical education due to simplicity of the
material model, lucidity of the basic ideas, and classical
character of mathematical machinery.

The paper is organized as follows.

In Section 2 the governing equations for von Mises ideal
rigid-plastic material under conditions of plain strain
deformation are introduced and a classic reduction to quasilinear
hyperbolic system of PDEs is given. Using {\it hodograph}
transformation, the system is transformed into the telegraph
equation, which means that  the original system is C-integrable.
This fact was known long ago \cite{kachanov}, but here the
transformation to the linear problem is done by appropriate change
of variables and does not involve any considerations from
mechanics. In Section 3 three main BVPs typical for applied
plasticity are defined. Using Riemann function for the telegraph
equation, the initial characteristic problem with analytic data
can be solved in terms of Bessel functions \cite{geiringer}. We
show that this solution can be expressed in terms of the
hypergeometric functions depending on the product of independent
variables,  which is more practical for computations. Exact
solutions for Cauchy and mixed problems are found by reducing the
problems to equivalent initial characteristic problem. In Section
4 the method is applied to two test problems from the theory of
plastic flow:
\begin{itemize}
\item Calculation of the slip-line field generated by two circular
arcs, and \item Computation of the stress state near an elliptic
hole loaded with constant normal pressure.
\end{itemize}
 Solutions obtained for the
test problems are in excellent agreement with those known from
literature.

To simplify exposition, the following notations and conventions
will be used hereafter.
\begin{itemize}

\item $t_n$= $t^n/n!,$ where $t$ is a variable.
 \item $[t]=[1,t_1,t_2,t_3,\ldots]^T$ (vector-column).
 \item Any given function is real analytic (i.e., it can be represented as
the sum of a power series). \item Summation is taken over the
range from $0$ to $\infty$ unless specified otherwise. \item
Notation
$$\,_0F_1\,(n+1;z)=n!\sum \frac{z^k}{k!(k+n)!}$$ stands for the
hypergeometric functions. The subindexes will be omitted for
simplicity. \item Given a power series $\sum a_k\,t_k,$ notation
$\mathbf{a}$ stands for the row of its coefficients.
 \end{itemize}

\section{Governing equations}

 Let a homogeneous isotropic body made of ideal rigid-plastic material be in a state
 of plane strain plastic deformation. Assume that the flow is parallel to $xy$-plane. Let $D$ denote
 the projection of the body onto this plane. It is assumed that the stress tensor components
 $\sigma_{ij},\ i,j=1,2,$
  and the strain rate tensor components $\epsilon_{ij}$ in the domain $D$ satisfy the following equations:
  \begin{enumerate}
\item Equilibrium equations
\begin{equation}
\sigma_{ij,\,i}=0,\ j=1,2. \label{eq:equilib}
\end{equation}
\noindent Here the Einstein rule of summation over repeated
indices is assumed.
 \item von Mises yield criterion
\begin{equation}
\left(\sigma_{11}-\sigma_{22}\right)^2+4\sigma_{12}^2 =4k^2,
\label{eq:mises}
\end{equation}
    where $k$ is the shear yield stress of the material.
\item  Constitutive relations
\begin{equation}
\frac{\sigma_{11}-\sigma_{22}}{\sigma_{12}}=\frac{\epsilon_{11}-\epsilon_{22}}{\epsilon_{12}},\
j=1,2. \label{eq:state}
\end{equation}
\item Incompressibility constraint
\begin{equation}  \label{eq:incomp}
\epsilon_{11}+\epsilon_{22}=0.
\end{equation}
\end{enumerate}
Two perpendicular lines, called slip-lines, are passing through
any interior point of the plastic domain. Each of these lines is
tangent to a maximum shear stress direction at this point. Lines
in {\it the first and the second shear directions} [7] are called
$\alpha$- and $\beta$-lines, respectively.  Under the plane strain
conditions, the mean stress $\sigma$ is given by
$\sigma=(\sigma_{11}+\sigma_{22})/2$. The stress components at a
point $(x,y)$ are determined by the equations
\begin{equation}
\sigma_{11}=\sigma-k\sin 2\phi,\ \sigma_{22}=\sigma+k\sin 2\phi,\
\sigma_{12}=k\cos 2\phi \label{eq:stress_tensor},
\end{equation}
where $\phi=\phi(x,y)$ is the angle from the $x$-axis to the
$\alpha$-line passing through the point.

Substituting these equations into (\ref{eq:equilib}) yields
 \begin{eqnarray} \label{eq:original_sys1}
\sigma_x - 2k\left(\phi_x\cos 2\phi+\phi_y\sin 2\phi\right)&=&0,\\
\label{eq:original_sys2} \sigma_y - 2k\left(\phi_x\sin
2\phi-\phi_y\cos 2\phi\right)&=&0, \label{eq:quasilin}
\end{eqnarray}
where lower indices indicate corresponding partial derivative.
Quasilinear system (3) is hyperbolic, and its characteristics
coincide with slip-lines. It is shown below that the system can be
converted into a linear PDE by a hodograph transform reversing the
roles of the dependent, ($\sigma,\ \phi$), and independent, ($x,\
y$), variables.

 Assuming that the Jacobian
 $$J=\sigma_x \phi_y-\sigma_y \phi_y$$
 does not vanish,
 one derives
 $$\phi_y=J\,x_{\sigma},\qquad \phi_x=-J\,y_{\sigma},\qquad
 \sigma_x=J\,y_{\phi}\qquad \sigma_y=-J\,x_{\phi},$$
and the system (\ref{eq:quasilin}) transforms into
 \begin{eqnarray*}
x_{\phi}&=& 2k\left( x_{\sigma}\cos 2\phi +y_{\sigma}\sin
2\phi\right),\\
y_{\phi}&=&2k\left(x_{\sigma}\sin 2\phi -y_{\sigma}\cos
2\phi\right).
\end{eqnarray*}

Introduce point-dependent rectangular coordinate system with axes
directed along characteristics,
\begin{equation}  \label{eq:moving_sys}
X=x\cos \phi +y\sin \phi,\qquad Y=-x\sin \phi +y\cos \phi,
\end{equation}
 and the characteristic coordinates $(\alpha,\beta),$
\begin{equation}  \label{eq:char_coord}
\alpha=\phi/2+(\sigma-\sigma_0)/(4k),\qquad
\beta=\phi/2-(\sigma-\sigma_0)/(4k),
\end{equation}
where $\sigma_0$ is the value of the mean stress at the origin.
Then, after algebraic simplifications, the governing system of
equations takes the form
\begin{equation}  \label{eq:pre_telegraphy}
Y_{\alpha}+X=0,\qquad X_{\beta}-Y=0.
\end{equation}

Geometric meaning of the characteristic coordinates   is clear
from Fig.1. If the functions $X,\ Y$ are smooth enough, each of
them satisfies the telegraphy equation
 \begin{equation}  \label{eq:telegraphy}
\frac{\partial f}{\partial\alpha\partial\beta}+f=0.
\end{equation}

\begin{figure}
\includegraphics*[scale=.4]{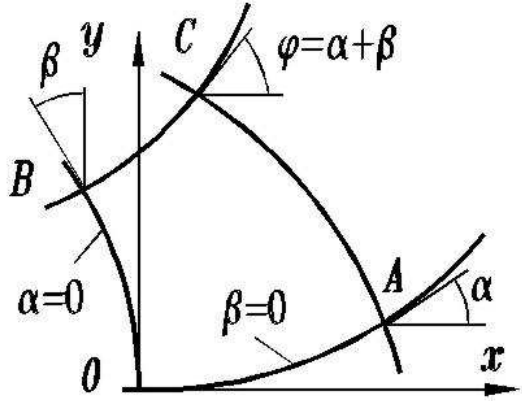}
\caption{ Curvilinear characteristics rectangle with initial
$\alpha$-line $OA$ and initial $\beta$-line $OB.$}
\end{figure}

The same  equation holds for the radii $R,\ S$ of  $\alpha$- and
$\beta$-line curvatures, respectively,  and for the components of
the velocity field \cite{sokol}. The stress components at a point
$(x(\alpha,\beta),y(\alpha,\beta))$ are defined by
(\ref{eq:stress_tensor}), with
 \begin{equation}  \label{eq:angle_stress}
\phi=\alpha+\beta,\qquad \sigma=\sigma_0+2k(\alpha-\beta).
\end{equation}
The main stresses are given by the equations $\sigma_1=\sigma+k,\
\sigma_2=\sigma-k$.
\section{Main boundary value problems}
\subsection{The initial characteristic problem}

It is well known  that there exists a unique solution to the
following initial characteristics problem \cite{courant}:

{\it Find a function $f(\alpha,\beta)$ satisfying
(\ref{eq:telegraphy}) in the domain $OACB$ (Fig.1) and the initial
conditions}
$$f(\alpha,0)=\sum c_n\alpha_n,\ \alpha\in(0,\alpha_1),\qquad
f(0,\beta)=\sum d_n\beta_n,\ \beta\in(0,\beta_1).$$ The
coefficients $c_n,\ d_n$ are given real numbers and the
compatibility condition $c_0=d_0$ holds. A classic solution in
terms of Bessel functions \cite{geiringer} can be rewritten in
terms of hypergeometric functions as
 \begin{equation}  \label{eq:sol_hyper}
f(\alpha,\beta)=\sum
\left(c_n\alpha_n+d_n\beta_n\right)\,F\,(n+1;-\alpha\beta)-c_0\,F(1;-\alpha\beta).
\end{equation}

 For a plasticity problem with curvatures of the initial characteristics
  given by the equations
 \begin{equation}  \label{eq:ini_radii}
R(\alpha,0)=\sum a_n\alpha_n,\qquad S(0,\beta)=\sum b_n\beta_n,
\end{equation}
missing  data $R(0, \beta),\ S(\alpha,0)$ are obtained from
equations (\ref{eq:pre_telegraphy}) (with $X,\ Y$ replaced with
$R,\ S$)
 \begin{equation}  \label{eq:missing_radii}
R(0,\beta)=\sum b_n\beta_{n+1}+a_0,\ S(\alpha,0)=\sum
a_n\alpha_{n+1}+b_0.
\end{equation}
Then it follows from (\ref{eq:sol_hyper}) that the slip-line field
in the characteristic rectangle $OACB$ is defined by the
curvatures
 \begin{eqnarray}  \label{eq:hyper_radius_r}
R(\alpha,\beta)=\sum\left(
a_n\alpha_n+b_{n-1}\beta_n\right)\,F(n+1;-\alpha\beta),\\\label{eq:hyper_radius_s}
S(\alpha,\beta)=\sum \left(-a_{n-1}\alpha_n+b_{n}\beta_n\right)\,
F(n+1;-\alpha\beta),
\end{eqnarray}
where  $a_{-1}=0=b_{-1}=0$.

The following theorem summarizes these results.

\newtheorem{riemann}{Theorem}

\begin{riemann} \label{riemann:box}
Let two arcs of intersecting slip-lines be given by the equations
(\ref{eq:ini_radii}). Then the slip-line field in the
characteristic rectangle generated by the two slip-lines is
uniquely defined by series (\ref{eq:hyper_radius_r}),
(\ref{eq:hyper_radius_s}).

Furthermore, if the coefficients $a_k,\ b_k$ are bounded, the
series  converge as exponential series. If for some $q>0$
$$|a_n|\le q^n\,n!,\qquad |b_n|\le q^n\,n!,$$
then for any $p,\ 0<p<1/q,$ the series converge  as geometric
series with ratio $r=pq$ provided $|\alpha|\le p,\ |\beta|\le p.$
\end{riemann}

 The estimates of the rate of convergence follow from the
 inequality
 $$ |F(n+1;z)|\le \exp(|z|).$$

{\it Remark.} A particular slip-line field generated by two
circular arcs is often used in technological plasticity. In this
case, $R(\alpha,0)=a_0,\ S(0, \beta)=b_0$ and solution takes the
form
$$R=a_0\,F\,(1;-\alpha\beta)+b_0\beta\,F\,(2;-\alpha\beta),\
S=b_0\,F\,(1;-\alpha\beta)-a_0\alpha\,F\,(2;-\alpha\beta).$$

Solution (\ref{eq:hyper_radius_r}), (\ref{eq:hyper_radius_s}) to
the initial characteristic problem can be written as  double power
series
\begin{eqnarray}  \label{eq:series_radii}
R(\alpha,\beta)&=&\sum\limits_{n,k}a_n\alpha_{n+k}(-\beta)_k+b_n\beta_{n+k+1}(-\alpha)_k,\\
S(\alpha,\beta)&=&\sum\limits_{n,k}-a_n\alpha_{n+k+1}(-\beta)_k+b_n\beta_{n+k}(-\alpha)_k,
\end{eqnarray}
or in a matrix
 form as
 \begin{equation}  \label{eq:matrix_radii}
R(\alpha,\beta)=\left(\mathbf{a}A(\beta)+\mathbf{b}B(\beta)\right)
[ \alpha],\
S(\alpha,\beta)=\left(-\mathbf{a}B(\alpha)+\mathbf{b}A(\alpha)\right)
[ \beta].
\end{equation}
Here $\mathbf{a,\ b}$ are the  vectors of coefficients in the
initial conditions (\ref{eq:ini_radii})
 and the matrix-functions $A,\ B$ are defined by the
formulas
\[A(t)=\left[ \begin{array}{ccccc}
               1&- t & t_2 & -t_3 & \ldots\\
                0 & 1&- t & t_2&\ldots \\
                0&0&1&  -t &\ldots \\
                0  &0&0&1& \ldots\\
                \ldots&\ldots&\ldots&\ldots&\ldots
                \end{array}
                \right],   \]
\[B(t)=\left[ \begin{array}{ccccc}
               t_1&- t _2& t_3 & -t_3 & \ldots\\
                t_2 &-t_3& t_4 &- t_5&\ldots \\
                t_3&-t_4&t_5& - t_6 &\ldots \\
                t_4  &-t_5&t_6&-t_7& \ldots\\
                \ldots&\ldots&\ldots&\ldots&\ldots
                \end{array}
                \right].   \]

If an initial slip-line shrinks to a point, the characteristic
rectangle degenerates into a fan with slip-line field defined by
(\ref{eq:matrix_radii}) with appropriate row of coefficients set
to zero.

\subsection{The Cauchy problem}
\begin{figure}[tbp]
\includegraphics*[scale=.6]{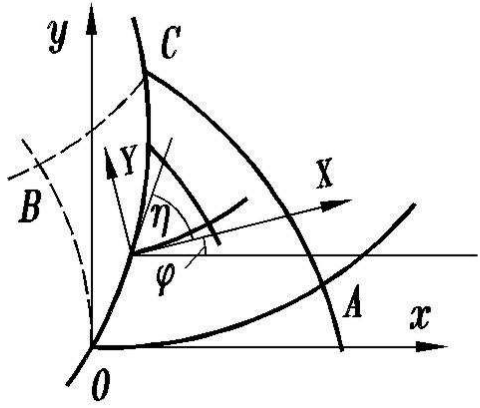}
\caption{ Characteristic triangle with Cauchy data on
non-characteristic curve $OC.$ }
\end{figure}
Consider a smooth non-characteristic curve $OC.$   Let the radius
of curvature of  $OC$  be $r(\gamma), \ \gamma$ being the angle
between the curve and $x-$axis. Given two functions
$\sigma(\gamma)$, $\phi(\gamma)$, the problem is to find functions
$\sigma(x,y)$, $\phi(x,y)$ satisfying equations
(\ref{eq:original_sys1}), (\ref{eq:original_sys2}) in the
 characteristic triangle $OAC$ (Fig.2) and boundary
conditions $\sigma=\sigma(\gamma)$, $\phi=\phi(\gamma)$ on  $OC$.
There exists a unique solution of this problem in the curvilinear
triangle $OAC$ \cite{courant}.

This BVP can be reduced to the equivalent initial characteristic
problem on the fictitious characteristic rectangle $OACB.$
Coefficients $a_k,\ b_k$  for the initial slip-lines are found as
follows. From an infinitesimal characteristic triangle with
hypotenuse lying on $OC$ (Fig. 2) one has the following elementary
identities
\begin{equation} \label{eq:rad along}
ds_{\alpha}/ \cos(\eta)=ds_{\beta}/ \sin(\eta)=r(\gamma)d\gamma,
\end{equation}
\noindent where  $ds_{\alpha},\ ds_{\beta}$ are the length
differentials of the bounding $\alpha$- and $\beta$-lines, and
$\eta=\gamma-\varphi$. Using (4), functions
$\alpha=\alpha(\gamma)$ and $\beta=\beta(\gamma)$ can be
determined at all points of $OC$. Substituting the expressions
\begin{equation}
ds_{\alpha}=R(\alpha(\gamma),\beta(\gamma))d\alpha,\qquad
ds_{\beta}=-S(\alpha(\gamma),\beta(\gamma))d\beta
\end{equation}
with $R, S$ defined by (\ref{eq:hyper_radius_r}),
(\ref{eq:hyper_radius_s}) into (\ref{eq:rad along}), one obtains a
system for the unknown coefficients $\alpha_k,\ \beta_k$. A
particular case when both the normal stress $\sigma_n$ and the
tangential stress $\tau_n$ are constant along $OC$ is detailed
below. In this case
\begin{equation}
\alpha=\beta, \ \eta=\pi/2-0.5 \cos^{-1}(\tau_n/k)=const,\
\gamma=\eta+2\alpha.
\end{equation}
Writing the radius of curvature of  $OC$  as a function of
$\alpha,$
$$r(2\alpha+\eta)=
\mathbf{r}[\alpha],$$ one obtains from (\ref{eq:rad along})
\begin{equation}\label{eq:rands}
R(\alpha,\alpha)/ \cos\eta=-S(\alpha,\alpha)/ \sin\eta
=\pm2\mathbf{r}[\alpha].
\end{equation}
These identities (with $R,\ S$ determined by
(\ref{eq:hyper_radius_r})), (\ref{eq:hyper_radius_s}))
  imply after
some mathematical manipulations
\begin{equation}\label{eq:a and b}
\mathbf{a}=2\mathbf{r}C, \qquad \mathbf{b}=2\mathbf{r}D,
\end{equation}
where the upper triangular matrices $C,\ D$, derived from the
recurrent relations in \cite{novo}, read
\[C=\left[ \begin{array}{rrrrrrr}
               c & s & c & s & c & s & \ldots\\
               0 & c & s & 2c & 2s & 3c & \ldots \\
               0 & 0 & 2!c & 1\cdot 2s & 2\cdot 3c & 2\cdot 3 s  &\ldots \\
               0 & 0 & 0 & 3!c & 1\cdot 2\cdot 3\cdot s& 2\cdot 3\cdot 4\cdot c & \ldots\\
               0 & 0 & 0 & 0 & 4!c & 1\cdot 2\cdot 3\cdot 4\cdot s
               & \ldots \\
                0 & 0 & 0 & 0 & 0 & 5!c & \ldots \\
                \ldots&\ldots&\ldots&\ldots&\ldots &\ldots&\ldots
                \end{array}
                \right],   \]
\[D=\left[ \begin{array}{rrrrrrr}
               -s & c & -s & c & -s & c & \ldots\\
               0 & -s & c & -2s & 2c & -3s & \ldots \\
               0 & 0 & -2!s & 1\cdot 2c & -2\cdot 3s & 2\cdot 3 c  &\ldots \\
               0 & 0 & 0 & -3!s & 1\cdot 2\cdot 3\cdot c& -2\cdot 3\cdot 4\cdot s & \ldots\\
               0 & 0 & 0 & 0 & -4!s & 1\cdot 2\cdot 3\cdot 4\cdot c
               & \ldots \\
                0 & 0 & 0 & 0 & 0 & -5!s & \ldots \\
                \ldots&\ldots&\ldots&\ldots&\ldots &\ldots&\ldots
                \end{array}
                \right].   \]
 For simplicity, notation c, s have been used for $\cos\eta$ and $\sin\eta$,
respectively.  Formulas (\ref{eq:a and b}) determine initial data
for the initial characteristic problem equivalent to the original
Cauchy problem in the domain $OAC$.

The solution for Cauchy problem is tested below on a classic
example with known slip-line field formed by logarithmic spirals.

\textbf{Example}. Consider a circular arc of radius one loaded
with constant normal pressure
  and zero tangent stress component.  This implies $\eta=\pi/4$  and $\mathbf{r}
=(1, \ 0, \ 0, \ ...).$ It follows from (\ref{eq:a and b}) that
$$\mathbf{a}=\sqrt{2}(1, \ 1, \
1, \ ...)\qquad \mathbf{b}=\sqrt{2}(1, \ -1, \ 1, \ -1, ...).$$
 Therefore the radii of curvature for the initial
$\alpha$- and $\beta$-lines are $R(a,0)=\sqrt{2}\exp(\alpha)$ and
$S(0,\beta)= \sqrt{2}\exp(-\beta)$). These are the radii of
curvature of logarithmic spirals.

 Stress-free surface
boundary $(\sigma_n=\tau_n=0$ on $OC$) can be analyzed similarly.
However, in this case the shape of the free boundary is not known
a priori, and to start a process of solution based on the
identities (\ref{eq:rands}) one needs to know one of the initial
slip-lines. Assuming for certainty that the initial slip-line $OA$
(Fig. 2) is known (for example, vector $\mathbf{a}$ could be found
from solution of BVP on the adjacent patch), one finds the row
$\mathbf{b}$ of coefficients for the initial $\beta$-line and free
surface radius of curvature $r(\alpha)$ from the identities
(\ref{eq:rands}). It can be shown that
\begin{equation}\label{eq:b and r}
\mathbf{b}=\mathbf{a}F, \qquad
r(2\alpha+\eta)=R(\alpha,\alpha)/\sqrt{2},
\end{equation}
where  the entries of the matrix $F$ are defined as
\begin{equation}\label{eq:free_surf}
f_{ii}=-1, \ f_{ij}=0,\ for j<i, \ f_{ij}/2=(-1)^{i+j+1} \ for \
j>i,
\end{equation}
\noindent and $R(\alpha,\beta)$ is given by
(\ref{eq:hyper_radius_r}). Finding the shape of a free surface may
be of interest in applications.

\subsection{Mixed problem}
Consider a smooth non-characteristic curve $OC$.  Let the radius
of curvature of  $OC$  be $r(\gamma), \ \gamma$ being the angle
between the curve and $x-$axis. Given the radius of curvature for
the initial $\alpha$-line $OA$ and a function $\phi(\gamma)$ on
$OC$, the problem is to find functions $\sigma(x,y)$, $\phi(x,y)$
satisfying equations (\ref{eq:original_sys1}),
(\ref{eq:original_sys2}) in the  characteristic triangle $OAC$ and
boundary condition $\phi=\phi(\gamma)$ on $OC$ (Fig.3).  It is
known that the problem has a unique solution \cite{courant}.

\begin{figure}[tbp]
\includegraphics*[scale=.4]{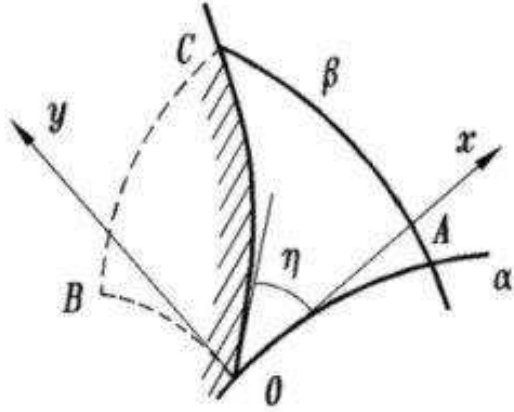}
\caption{Characteristic triangle for the mixed problem. One of the
bounding slip-lines and the values of the function $\varphi$  on
the non-characteristic curve $OC$ are given.}
\end{figure}
Such formulation arises in applied plasticity if a tangential
stress $\tau_n$ is given on a contact line $OC$. It is assumed
here that $\tau_n=\mu k, \ \mu\in[0,1],$ (Prandtl friction law).
Note that constancy of $\tau_n$ implies that the angle $\eta$ also
remains constant on $OC.$ Assuming that the curve $OC$ is defined
by an equation $\beta=\beta(\alpha)=\mathbf{c}[\alpha],$
 one obtains elementary identities
\begin{equation}\label{eq:thirteen}
R(\alpha,\beta(\alpha))/\cos\eta=-S(\alpha,\beta(\alpha))\beta/\sin\eta=r(\alpha+\beta(\alpha)+\eta)(1+\beta').
\end{equation}
The unknown vector $\mathbf{c}$  and vector $\mathbf{b}$ for the
initial $\beta$-line of the equivalent
 initial characteristic problem can be obtained from these
identities. A particular case  is analyzed below.

Let  $OC$ be a straight line. If $x$-axis is directed along the
$\alpha$-line passing through the origin,  then
$\varphi=\alpha+\beta=0$
 on $OC$,  or equivalently, $\beta=-
\alpha$. The rightmost term in (\ref{eq:thirteen}) takes the form
$\infty\cdot 0$ and becomes useless. The first identity takes the
form
\begin{equation}\
R(\alpha,-\alpha)/\cos\eta=-S(\alpha,-\alpha))/\sin\eta
\end{equation}
Using (\ref{eq:series_radii}),  one obtains after some
mathematical manipulations the missing boundary data for the
initial $\beta$-line of the equivalent initial characteristic
problem. Specifically, if $ \eta<\pi /2$, the row of the
coefficients for the $\beta$-line radius of curvature is given by
the equation $b=aT(\eta)$ , where the elements of the upper
triangular matrix $T$ are defined by the equations
\begin{eqnarray*}\
t_{ij}&=&(-1)^i(\tan\eta)^{j-i-1}(\tan^2\eta-1) \ for\ j>i, \\
t_{ii}&=&(-1)^i\tan\eta \ for \ i=j\\
t_{ij}&=&0\ for\ j<i.
\end{eqnarray*}
 It can be shown that for perfectly rough boundary $(\eta
 =\pi/2)\   b_n=(-1)^n a_n.$

\section{Test problems}
\begin{figure}[tbp]
\includegraphics*[scale=.4]{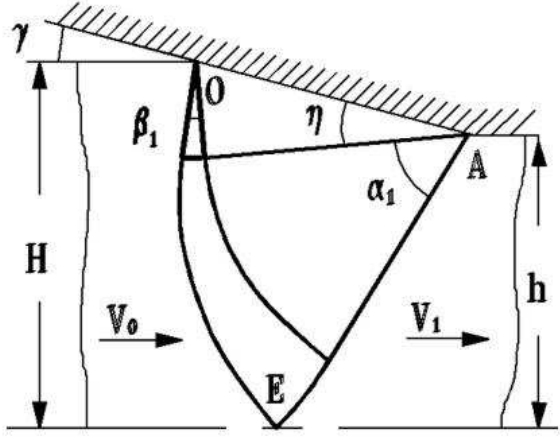}
\caption{Slip-line field of extrusion/drawing through the short
wedge-shaped die.}
\end{figure}

In this section SLAT is applied to two test problems. The first
example is classical and its exact solution is well known
\cite{nepershin}. The second problem was treated before only
numerically. Since the main mathematical operation used in SLAT is
multiplication of  matrices and vectors of small size, computer
time needed for each problem is negligible. However, some
preliminary analytic work is needed to present initial data in the
required form.

\subsection{Extrusion through a short wedge-shaped die}

Consider an extrusion   of a plastic material through a short
wedge-shaped die of angle $\gamma$ (Fig. 4). For comparability of
this example with solution given in \cite{johnson}, the die is
assumed to be frictionless (i.e., $\eta=\pi/4$), and the following
values of parameters are used: $\gamma=10^o,\ \alpha_1=30^o,\
\beta_1=-20^o,\ OA=2$ (dimensionless units). Thus the initial data
are $\mathbf{a}=(\sqrt{2}, \ 0, \ 0,...)$ and
$\mathbf{b}=(-\sqrt{2}, \ 0, \ 0,...)$. The initial characteristic
problem with the data
 was solved in the characteristic rectangle (Fig. 4)
using (\ref{eq:ini_radii}). Then the coordinates of the point $E$
relative to the origin $O$ were found using simple geometric
considerations. The hydrostatic stress components $p_B=-\sigma(B),
\ p_D=-\sigma(D)$, and the extrusion pressure $P/H$ were obtained
using the solution for the initial characteristic problem and
Hencky`s first theorem  in a  standard static manner
\cite{johnson}. The value of $H$ in this example is found to be
$2.28774.$ \vspace{0.5cm}

 Table 1.Comparison of some parameters for the extrusion process

\begin{tabular}{c|c|c|c}
    \hline\\
  & $\mathbf{x_E}$ & $\mathbf{y_E}$ & $\mathbf{P/H}$\\
  \hline\\
  \cite{johnson} & $0.9065$ &  $-2.2877$ & $ 0.4117$ \\
  \hline \\
  $SLAT$ & $0.90648$ & $-2.28774$ &    $0.41164$ \\
  \hline
  \end{tabular}

\vspace{0.5cm}

For accuracy of $10^{-5}$, five-dimensional truncation of matrices
and vectors was found to be sufficient. The results for some
parameters obtained by the original matrix method \cite{johnson}
and SLAT are practically identical (Table 1). Hydrostatic and
extrusion pressures are normalized by the shear yield stress of
the material.

 \subsection{Stress state calculation near loaded elliptic
 hole}

 Consider planar plastic flow near elliptic hole loaded with constant normal pressure (Fig. 5).
 Firstly,  Cauchy problem in the domain 1 was solved.
 Secondly, initial characteristic problem was solved in the domain 2 (note that the domain 2 is the upper half
 of a characteristic rectangle).
\begin{figure}[tbp]
\includegraphics*[scale=.5]{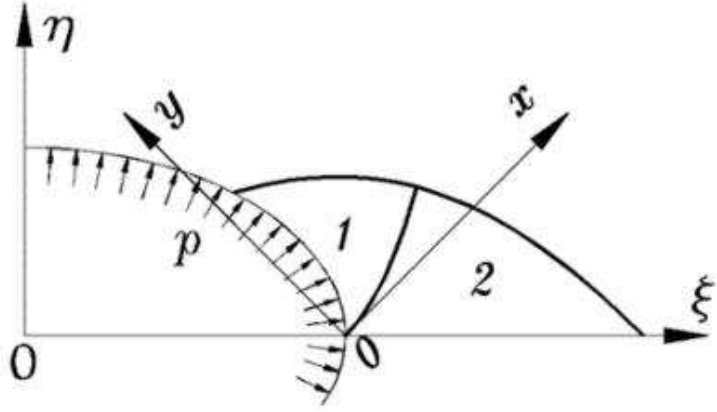}
\caption{Slip-line field for one quarter of an elliptic hole
loaded with constant normal pressure.}
\end{figure}
Equation for the radius of curvature of the elliptic contour is
obtained as follows. Let lengths be normalized by the semimajor
axis of the ellipse. Then using parametric equations $x=\cos t,\
y=b\sin t,\ 0<b\le 1,$ of the ellipse and the standard formula for
the curvature, one obtains  the radius of curvature $\rho (t)$ in
the form
\begin{equation}\
\rho (t)=2\sqrt{2}((1+b^2)+(1-b^2)\cos{2t})^{3/2}/b.
\end{equation}
From elementary geometry it follows that
\begin{equation}
\gamma =2\alpha+\pi/4=\tan^{-1}(b\cot t)-\pi/4.
\end{equation}
From this equation the variable t can be expressed in terms of
$\alpha$, and after some mathematical manipulations one  obtains
the radius of the curvature for the elliptic hole as a function of
$\alpha:$
\begin{equation}\
r(\alpha)=\rho(t(\alpha))=2\sqrt{2}b^2(1+d^2)^{-3/2}(1+q\cos4\alpha)^{-3/2},
\end{equation}
where $q=(1-b^2)/(1+b^2)$. For comparability with \cite{sokol},
value $b=0.4$ have been chosen. For the domain 1, the row of the
coefficients for the initial $\alpha$-line of the equivalent
initial characteristic problem was found using the first equation
in (\ref{eq:a and b}). For accuracy of $10^{-4}$,
fifteen-dimensional truncation of  vectors was found to be
sufficient. In general, numerical experimentations show that the
number of terms needed for a given accuracy grows significantly
when the value of $b$ decreases. Coefficients for the initial
$\beta$-line, which is symmetric to the $\alpha$-line with respect
to the horizontal axis, are given by the equation
$b_k=(-1)^{k+1}a_k.$ The main stresses on the $\xi$-axis
normalized by $2k$ are given by the equations
$\sigma_1=\sigma+(2\alpha+1), \ \sigma_2=\sigma+2\alpha$ . The
comparative results are shown in Fig. 6, where $\Delta=\sigma_1
+p.$

\begin{figure}[tbp]
\includegraphics*[scale=.5]{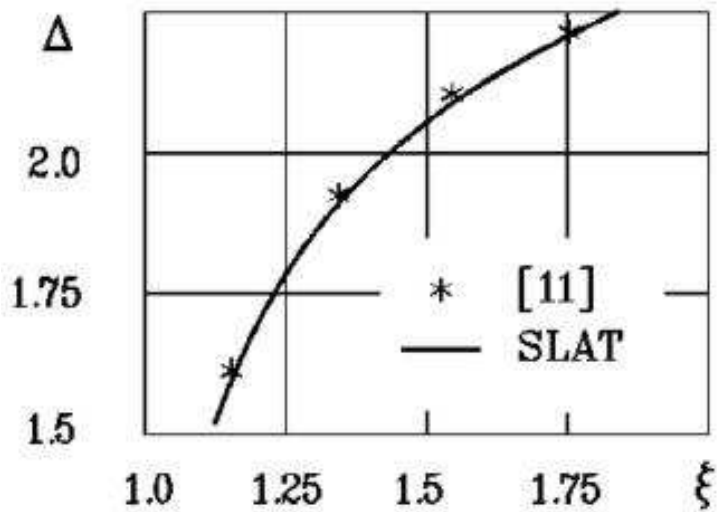}
\caption{Comparison of the computation of $\Delta=\sigma_1 +p$ for
the elliptic hole done numerically \cite{sokol} and by SLAT.}
\end{figure}
\section{Conclusions}

Mathematically accurate  and efficient analytic implementation of
the slip-line method aimed at computing the stress fields for
plane strain deformation of the rigid-plastic medium has been
presented. It can be used as a source of reliable test problems
for numerical methods and in  engineering applications.

 The velocity field in
the hodograph plane can also be found by this method. Using
solutions in the physical and hodograph planes, exact expressions
for the plastic power can also be derived. This gives an
alternative way for computing the technological pressure, an
important parameter for engineering applications, and also
provides a framework for solving relevant optimization problems.
Detailed explicit formulas for  computation of energetic
characteristics of slip-line fields will be  given elsewhere.

Although only Prandtl friction law was analyzed in  the mixed
problem presented in this work, Coulomb friction requires just
technical modifications and can be treated  in a similar way. This
may be of interest in the theory of granular flow, where using
this kind of friction law is customary \cite{drescher}. This
modification of SLAT will be presented in a later paper.
\subsection*{Acknowledgements}
This research was partially supported by the CSU/AAUP grant \#
242414.

\end{document}